\documentclass{amsart}

\usepackage{latexsym}
\usepackage{amsmath}
\usepackage{amssymb}
\usepackage{amsthm}
\usepackage[all]{xy}

\allowdisplaybreaks                          
\newtheorem{thm}{Theorem}
\newtheorem{prop}[thm]{Proposition}

\newtheorem{rem}[thm]{Remark}

\newtheorem{ex}[thm]{Example}

\begin{document}

\title{Normal projective varieties of degree $5$}

\author{Andrea Luigi Tironi}

\date{\today}

\maketitle

\begin{abstract}
We list all the irreducible reduced and not degenerate normal
projective varieties $X\subset\mathbb{P}_k^N$ of dimension $n$ and
degree five
defined over an algebraically closed field $k$ of
$\mathrm{char}(k)=0$. In the smooth case, or when $n=2$, we give
also classification results for any algebraically closed field
$k$ of $\mathrm{char}(k)\geq 0$.
\end{abstract}

\section{Introduction}

Let $X$ be an irreducible closed subvariety of degree $d$ of the
projective space $\mathbb{P}_k^N$ defined over an algebraically
closed field $k$.

If $d\leq 4$, then the list of all such varieties is known (see  
\cite{***} for $d=3$ and H.P.F. Swinnerton-Dyer \cite{S} for
$d=4$).

Let now $d=5$. If $X$ is a smooth complex projective variety, then
the classification is also known (see P. Ionescu \cite{ionescu}).
Furthermore, the classification of all complex normal surfaces was
obtained by E. Halanay in \cite{H1}. Finally, we also mention here that the classification of all
complex surfaces of degree $5$ in $\mathbb{P}_{\mathbb{C}}^N$
which are not contained in any hyperplane of $\mathbb{P}_{\mathbb{C}}^N$ is
known too (see, for example, G. Castelnuovo \cite{C1} and
\cite{C2}, J.G. Semple and L. Roth \cite[p. 218]{SR}, M. Nagata
\cite{N}, A. Lanteri and M. Palleschi \cite{LP} and C. Okonek
\cite{O}).

\smallskip

The purpose of this note is to show that the results of \cite{H1}
and \cite{ionescu2} hold also in higher dimensions and in positive
characteristic respectively. In doing so, we first use slight
different techniques than those in \cite{H1}, 
obtaining a list of all varieties $X\subset\mathbb{P}_k^N$ of dimensions $n\geq 2$ and degree $5$
when essentially either $\mathrm{char}(k)=0$ and $X$ is normal, or
$\mathrm{char}(k)\geq 0$ and $X$ is smooth or $n=2$.

\medskip

More precisely, in Section $3$ we prove the following

\medskip

\noindent\textbf{Theorem.} \textit{Let $X\subset\mathbb{P}_k^N$ be
an irreducible and reduced 
normal projective variety of 
dimension $n\geq 1$ defined over an algebraically closed field $k$
of $\mathrm{char}(k)=0$. 
Denote by $H$ 
an ample and spanned line bundle on $X$. If $d:=H^n=5$, then
either $X$ is a covering $\pi:X\to\mathbb{P}_k^n$ of degree five
over $\mathbb{P}_k^n$ with
$H=\pi^*(\mathcal{O}_{\mathbb{P}_k^n}(1))$, or $H$ is very ample
and one of the following cases can occur:
\begin{enumerate}
\item $n=2$ and $(X,H)$ is a scroll over an elliptic curve; \item
$n\leq 5$ and $(X,H)$ is a smooth scroll over $\mathbb{P}_k^1$,
i.e. $X\cong\mathbb{P}(\mathcal{O}(a_1)\oplus
...\oplus\mathcal{O}(a_n))$ and $H=\xi$ is the
tautological line bundle, where the $a_i$'s 
are positive integers such that $a_1+...+a_n=5$; \item $n\leq 5$
and $(X,H)$ is a normal Del Pezzo variety with
$K_X+(n-1)H\simeq\mathcal{O}_X$; \item $n=6$ and $X$ is the
grassmannian manifold $\mathbb{G}(1,4)\subset\mathbb{P}_k^9$ of
lines in $\mathbb{P}_k^4$ embedded by the Pl\"{u}cker embedding;
\item $X\subset\mathbb{P}_k^{n+2}$ is linked to a linear space
$P\cong\mathbb{P}_k^n$ by a complete intersection of type $(2,3)$,
i.e., $X\cup P=V_2\cap V_3$, where $V_3\subset\mathbb{P}_k^{n+2}$
is a cubic hypersurface and $V_2\subset\mathbb{P}_k^{n+2}$ is a
quadric cone of rank $r$
with $3\leq r\leq 4$; 
moreover, $X$ is a Weil divisor on $V_2$ and one of the following possibilities can occurs:
\begin{enumerate}
\item[$(i)$] $r=3$, $X\simeq 5P_1$ and $\mathrm{Cl}(V_2)\cong\mathbb{Z}[P_1]$;
\item[$(ii)$] $r=4$, $X\simeq 2P_1+3P_2$ and $\mathrm{Cl}(V_2)\cong\mathbb{Z}[P_1]\oplus\mathbb{Z}[P_2]$,
\end{enumerate}
where $\mathrm{Cl}(V_2)$ is the class group of $V_2$ and the
$P_i$'s are linear spaces $\mathbb{P}_k^n$ such that $P_1\simeq
P$; furthermore, if $X$ is smooth then $n\leq 3$, and 
when $n=3$ we have only the case $(ii)$; \item $X$ is a
quintic hypersurface of $\mathbb{P}_k^{n+1}$; \item $(X,H)$ is a
generalized cone over a variety $V\subset X$ as in cases $(1)$ to
$(4)$.
\end{enumerate}}

\medskip

Finally, in Section $4$ we deduce from the proof of the above
Theorem some results which extend those in \cite{H1} and
\cite{ionescu2} to any algebraically closed field $k$ of
$\mathrm{char}(k)\geq 0$ (see Prop. 2 and Thm. 5).

\section{Preliminaries}

Let $X$ be a projective normal subvariety of dimension $n\geq 2$
of a projective space $P\simeq\mathbb{P}_k^N$ defined over an
algebraically closed field $k$ of characteristic zero such that
$H^0(P,\mathcal{O}(1))\to H^0(X,\mathcal{O}_X(1))$ is bijective.
Denote by $H$ the restriction of $\mathcal{O}_P(1)$ to $X$.
Without loss of generality, we can suppose that $X$ is a not
degenerate irreducible and reduced normal variety of $P$.

For basic definitions of the invariants degree $H^n$, sectional
genus $g(X,H)$ and $\Delta$-genus $\Delta (X,H)$ of a polarized
variety $(X,H)$ as above, we refer to \cite[$\S$I (2.0), (2.1),
(2.2)]{Fu1}.

Denote by $d:=\deg (X)=H^n$. Throughout this note, we assume that $d=5$.

Given two points $p,q$ on a projective space $P:=\mathbb{P}_k^N$,
we denote by $p*q$ the line passing through $p$ and $q$. For
subvarieties $V,W$ of $P$, $V*W$ denotes the closure of the union
of all the lines $p*q$ with $p\in V$ and $q\in W$. We define
$x*x=x$ and $V*\emptyset =V$ by convention. So, if $p$ is a point
of $P$, $p*X$ is the cone over $X$ with vertex $p$.

A point $x$ on a subvariety $X$ of $P$ is called a vertex of $X$
if $x*X=X$. The set of vertices of $X$ is called the ridge of $X$,
and is denoted by $\mathrm{Rdg}(X)$. The variety $X$ is said to be
conical if $\mathrm{Rdg}(X)\neq\emptyset$.

By ``generalized cone'' we mean the following. Let $\mathcal{A}$
be an ample vector bundle on a curve $C$ and set
$\mathcal{E}=\mathcal{O}_C^{\oplus s}\oplus\mathcal{A}$,
$M=\mathbb{P}_C(\mathcal{E})$ and $L=\mathcal{O}_M(1)$. The
surjection $\mathcal{E}\to\mathcal{O}_C^{\oplus s}$ defines a
submanifold N of M such that
$N\cong\mathbb{P}(\mathcal{O}_C^{\oplus s})\cong
C\times\mathbb{P}_k^{s-1}$ with $L_N$ being the pull-back of
$\mathcal{O}(1)$ of $\mathbb{P}_k^{s-1}$. The linear system $|mL|$
has no base point for $m>>0$ and gives a birational morphism $\pi:
M\to W$ onto a normal variety W such that $L=\pi^*A$ for some
ample line bundle $A$ on $W$. Moreover $M-N\cong W-\pi (N)$ and
$\pi_N:N\to\pi (N)$ is the second projection of $N$ onto
$\mathbb{P}_k^{s-1}$. If $(X,H)$ is isomorphic to such a polarized
variety $(W,A)$, it is called a generalized cone over
$(\mathbb{P}(\mathcal{A}),\mathcal{O}(1))$.

\section{Proof of the Theorem}

First of all, note that $H^n=\deg\varphi(X)\cdot\deg\varphi$ is a
prime number, where $\varphi:X\to\mathbb{P}_k^N$ is the morphism
associated to the linear system $|H|$. Thus we have either (i)
$\deg\varphi(X)=1, \deg\varphi=5$, or (ii) $\deg\varphi(X)=5,
\deg\varphi=1$. In case (i), we see that $X$ is a covering
$\pi:X\to\mathbb{P}_k^n$ of degree five over $\mathbb{P}_k^n$ with
$H=\pi^*(\mathcal{O}_{\mathbb{P}_k^n}(1))$, while in case (ii) $H$
is very ample.

If $X$ is conical, then by \cite[$\S$I (5.13)]{Fu1} it suffices to
consider $X\cap L$ for a suitable linear subspace $L$ of
$\mathbb{P}_k^{N}$ because $X=R*(X\cap L)$, where
$R=\mathrm{Rdg}(X)$. So, from now on, we can assume without loss
of generality that $X$ is not conical.

From the $\Delta$-genus formula it follows that
\begin{equation}
4=\Delta + N - n, \tag{*}
\end{equation}
where $\Delta :=n+H^n-h^0(H)$ is the $\Delta$-genus of the
polarized pair $(X,H)$ and $h^0(H)=N+1$. Note that $\Delta\geq 0$
and $N - n\geq 1$. In particular, we have $\Delta\leq 3$ with
equality if and only if $X\subset\mathbb{P}_k^{n+1}$ is a
hypersurface of degree 5. So, from now on, we can assume that
$\Delta\leq 2$.

\smallskip

First of all, note that $X$ is smooth if $n=1$. Thus by \cite{Fu1} and \cite{ionescu},
if $n=1$ we deduce that $X$ is one of the following curves:
\begin{enumerate}
\item[(a)] a plane curve of genus $6$; \item[(b)] a curve of genus
$2$ in $\mathbb{P}_k^3$; \item[(c)] the intersection of
$\mathbb{G}(1,4)\subset\mathbb{P}_k^9$ with five general
hyperplanes;
\item[(d)] the rational normal curve in $\mathbb{P}_k^5$. 
\end{enumerate}

\smallskip

\noindent From now on suppose that $n\geq 2$.

\smallskip

Having in mind that $X$ is assumed to be not conical, if
$\Delta\leq 1$ then by \cite[$\S$I (5.10), (5.15), (6.9), (8.11),
$\S 9$]{Fu1} and \cite{ionescu} it follows that $(X,H)$ is one of
the following pairs:
\begin{enumerate}
\item[$(a')$] an scroll over $\mathbb{P}_k^1$ with $n\leq 5$, i.e.
$(X,H)\cong (\mathbb{P}(\mathcal{O}(\alpha_1)\oplus
...\oplus\mathcal{O}(\alpha_{n})),\xi)$, where $\xi$ is the
tautological line bundle and all the $\alpha_i$'s are positive
integers such that $\alpha_1+...+\alpha_{n}=5$; \item[$(b')$] a
linear section of $\mathbb{G}(1,4)\subset\mathbb{P}_k^9$;
\item[$(c')$] a normal Del Pezzo variety with $n\leq 5$.
\end{enumerate}

\smallskip

\noindent Thus $\Delta =2$ and by (*) we obtain that $N=n+2$.
Furthermore, since $X$ is not conical, if the sectional genus
$g:=g(X,H)$ of the pair $(X,H)$ is less than or equal to one, then
from \cite[$\S$I (3.4); $\S$III (18.21)]{Fu1} we know that $(X,H)$
is a scroll over an elliptic curve. Hence $X\cong\mathbb{P}_C(E)$
and $H=\mathcal{O}(1)$, where $E$ is an ample and spanned vector
bundle of rank $n$ over an elliptic curve $C$. Let $\pi :X\to C$
be the scroll map. By applying $\pi_*$ to the exact sequence
$$0\to\mathcal{O}_X\to\mathcal{O}_X(H)\to\mathcal{O}_H(H)\to 0,$$ we obtain by \cite[$\S$III Ex.8.4]{Ha}
\begin{equation}
0\to\mathcal{O}_C\to\pi_*\mathcal{O}_X(H)\to\pi_*\mathcal{O}_H(H)\to R^1\pi_*\mathcal{O}_X=0, \tag{**}
\end{equation}
where $\pi_*\mathcal{O}_X(H)=E$ is ample on $C$. Therefore we have
$H^1(C,E)\cong H^0(C,E^*)=0$, where $E^*$ is the dual of $E$. Thus
the cohomology sequence of (**) gives
$h^0(C,\pi_*\mathcal{O}_X(H))=h^0(C,\pi_*\mathcal{O}_H(H))$ and by
\cite[$\S$III Ex. 8.1]{Ha} we conclude that
$h^0(X,\mathcal{O}_X(H))=h^0(H,\mathcal{O}_H(H))$. Since $\Delta
=1+\deg H_C-h^0(H_C)=1$ for an elliptic curve $C$, we have by
induction $\Delta (H,\mathcal{O}_H(H))=n-1$, i.e.
$h^0(H,\mathcal{O}_H(H))=H^n$. Thus it follows that $\Delta
(X,\mathcal{O}_X(H))=n+H^n-h^0(X,\mathcal{O}_X(H))=n$, that is,
$n=\Delta =2$. 

Now, we can assume that
$$N=n+2,\ \ \Delta=2 \ \ \mathrm{and} \ \ g\geq 2.$$
Since $g\geq 2=\Delta$, by \cite[$\S$I (3.5)]{Fu1} we see that
$g=\Delta=2$. In particular, $X$ can be assumed to be not
contained in any hyperplane of $\mathbb{P}_k^{n+2}$, i.e. not
degenerate.

\smallskip

Take $n-i$ general hyperplanes $H_j$ and put $X_i:=X\cap H_1\cap
...\cap H_{n-i}\subset\mathbb{P}_k^{i+2}$ for $i=1,...,n-1,$ and
$X_n:=X$. Note that $X_1$ is a smooth curve, since $X_n$ is
normal.

For any $k=2,...,n$, consider now the following exact sequences
$$0\to\mathcal{O}_{X_k}\to\mathcal{O}_{X_k}(1)\to\mathcal{O}_{X_{k-1}}(1)\to 0,$$
$$0\to\mathcal{O}_{X_k}(1)\to\mathcal{O}_{X_k}(2)\to\mathcal{O}_{X_{k-1}}(2)\to 0$$
and
$$0\to\mathcal{O}_{X_k}(2)\to\mathcal{O}_{X_k}(3)\to\mathcal{O}_{X_{k-1}}(3)\to 0.$$
By an induction argument, we deduce the following inequalities
$$h^0(\mathcal{O}_{X_k}(1))\leq h^0(\mathcal{O}_{X_k})+h^0(\mathcal{O}_{X_{k-1}}(1))
=1+h^0(\mathcal{O}_{X_{k-1}}(1))\leq ...\leq (k-1)+h^0(\mathcal{O}_{X_1}(1)),$$
$$h^0(\mathcal{O}_{X_k}(2))\leq h^0(\mathcal{O}_{X_k}(1))+...+h^0(\mathcal{O}_{X_2}(1))+
h^0(\mathcal{O}_{X_1}(2))\leq \frac{k(k-1)}{2}+(k-1)\cdot h^0(\mathcal{O}_{X_1}(1))+h^0(\mathcal{O}_{X_1}(2))$$
and
$$h^0(\mathcal{O}_{X_k}(3))\leq h^0(\mathcal{O}_{X_k}(2))+h^0(\mathcal{O}_{X_{k-1}}(3))\leq ...
\leq \sum_{j=2}^{k}h^0(\mathcal{O}_{X_j}(2))+h^0(\mathcal{O}_{X_{1}}(3))\leq  $$
$$\leq \sum_{j=2}^{k}\left[ \frac{j(j-1)}{2}+(j-1)\cdot h^0(\mathcal{O}_{X_1}(1))+h^0(\mathcal{O}_{X_1}(2))
\right]+h^0(\mathcal{O}_{X_{1}}(3))=$$
$$= \frac{1}{2}\left[ 1+\sum_{j=2}^{k} j^2-1-\sum_{j=2}^{k} j\right]+\left[\sum_{j=2}^{k} (j-1)\right]\cdot
h^0(\mathcal{O}_{X_1}(1))+(k-1)\cdot h^0(\mathcal{O}_{X_1}(2))+h^0(\mathcal{O}_{X_{1}}(3)),$$
i.e.,
\begin{equation}
h^0(\mathcal{O}_{X_k}(2))\leq \frac{k(k-1)}{2}+(k-1)\cdot
h^0(\mathcal{O}_{X_1}(1))+h^0(\mathcal{O}_{X_1}(2))\tag{$\alpha$}
\end{equation}
\begin{equation}
h^0(\mathcal{O}_{X_k}(3))\leq
\frac{k}{6}(k^2-1)+\frac{k(k-1)}{2}\cdot
h^0(\mathcal{O}_{X_1}(1))+ (k-1)\cdot
h^0(\mathcal{O}_{X_1}(2))+h^0(\mathcal{O}_{X_{1}}(3))\tag{$\beta$}
\end{equation}
Since $g(X_1)=2, \deg\mathcal{O}_{X_1}(1)=5$ and
$h^1(\mathcal{O}_{X_1}(m))=0$ for $m\geq 1$, by the Riemann-Roch
Theorem we deduce that
$$h^0(\mathcal{O}_{X_1}(1))=4, \ h^0(\mathcal{O}_{X_1}(2))=9, \ h^0(\mathcal{O}_{X_1}(3))=14.$$
Moreover, if $\mathcal{I}$ denotes the ideal sheaf of $X$ in
$\mathbb{P}_k^{n+2}$, from the following exact sequences
$$0\to \mathcal{I}(2)\to\mathcal{O}_{\mathbb{P}_k^{n+2}}(2)\to\mathcal{O}_X(2)\to 0$$
and
$$0\to \mathcal{I}(3)\to\mathcal{O}_{\mathbb{P}_k^{n+2}}(3)\to\mathcal{O}_X(3)\to 0,$$
it follows from $(\alpha)$ and $(\beta)$ that
\begin{equation}
h^0(\mathcal{I}(2))\geq
h^0(\mathcal{O}_{\mathbb{P}_k^{n+2}}(2))-h^0(\mathcal{O}_X(2))\geq
1,\tag{$\diamondsuit$}
\end{equation}
\begin{equation}
h^0(\mathcal{I}(3))\geq
h^0(\mathcal{O}_{\mathbb{P}_k^{n+2}}(3))-h^0(\mathcal{O}_X(3))\geq
n+5.\tag{$\diamondsuit'$}
\end{equation}
So by $(\diamondsuit)$ we conclude that
$$\dim H^0(\mathbb{P}_k^{n+2},\mathcal{I}(2))\geq 1.$$
An element of that space is a form of degree 2, whose zero-set
will be an hypersurface $Q\subset\mathbb{P}_k^{n+2}$ of degree $2$
containing X. It must be irreducible (and reduced), because $X$ is
irreducible and not degenerate. Moreover, $X$ could not be
contained in two distinct irreducible quadric surfaces $Q,Q'$,
because then it would be contained in their intersection $Q\cap
Q'$ which is a variety of degree $4$, and that is impossible
because $\deg X = 5$. So we see that $X$ is contained in a unique
irreducible quadric hypersurface $Q$.

Furthermore, from ($\diamondsuit'$) we deduce also that
$$\dim H^0(\mathbb{P}_k^{n+2},\mathcal{I}(3))\geq n+5.$$
The cubic forms in here consisting of the quadratic form above
times a linear form, form a subspace of dimension $n+3$. Hence
there exists at least an irreducible cubic form in that space. So
$X$ is contained in an irreducible cubic hypersurface
$V_3\subset\mathbb{P}_k^{n+2}$. Then $X$ must be contained in the
complete intersection $Q\cap V_3$, and since $\deg X=5$ and $\deg
(Q\cap V_3)=6$, this shows that $X$ is linked to a linear space
$P:=\mathbb{P}_k^n$ by a complete intersection $Q\cap V_3$ of type
$(2,3)$ in $\mathbb{P}_k^{n+2}$.

Note that $Q$ must be singular, since $P\subset Q$ and $n\geq 2$.
Hence $Q$ is a quadric cone. After a suitable linear change of
variables, $P$ is defined by $x_0=x_1=0$ and $Q$ can be brought
into the form $x_0L_1+x_1L_2=0$, where the $L_i$'s are linear form
in the variables $x_0,...,x_{n+2}$. If the linear forms
$x_0,x_1,L_1,L_2$ are independent, then after a linear change of
variables, we can write the equation of $Q$ as (i) $x_0x_2+x_1x_3=0$. On the
other hand, if the linear forms $x_0,x_1,L_1,L_2$ are not
independent, then after a suitable linear change of variables, we
can write $Q$ into the form (ii) $x_0^2+x_1x_2=0$. Furthermore, in
case (ii), we have Cl$(Q)\cong\mathbb{Z}[P]$ and $H_Q\simeq 2P$. So we deduce that
$X+P\simeq 3H_Q\simeq 6P$, that is $X\simeq 5P$, where
$P\cong\mathbb{P}_k^n$ is the linear space on $Q$ which is linked
to $X$ by the complete intersection $Q\cap V_3$. Finally, 
in case (i) we have Cl$(Q)\cong\mathbb{Z}[P]\oplus\mathbb{Z}[P']$, where
$P'\cong\mathbb{P}_k^n$ is another linear space in $Q$ not
linearly equivalent to $P$. Therefore $H_Q\simeq P+P'$ and we can
conclude in this case that $X+P\simeq 3H_Q\simeq 3P+3P'$, that is
$X\simeq 2P+3P'$.

Assume now that $X$ is smooth. Since $V_3$ contains the linear
space $P: x_0=x_1=0$, we see that $V_3$ is described by an
equation of type $x_0F+x_1G=0$, where $F,G$ are two homogeneous
polynomials of degree two in the variables $x_0,...,x_{n+2}$.
Furthermore, note that $\mathrm{Sing}(V_3)\cap P$ is described by
$x_0=x_1=F=G=0$, i.e. a complete intersection $V_{2,2}$ of type
$(2,2)$ on $P\cong\mathbb{P}_k^n$, and that $\mathrm{Sing}(Q)\cap
P=\mathrm{Sing}(Q)$ is either a linear space
$\mathbb{P}_k^{n-2}\subset P$ with equations $x_0=x_1=x_2=x_3=0$
in case (i), or a linear space $\mathbb{P}_k^{n-1}\subset P$ with
equations $x_0=x_1=x_2=0$ in case (ii). Note that
$$\mathrm{Sing}(Q)\cap\mathrm{Sing}(V_3)\supseteq
(\mathrm{Sing}(Q)\cap P)\cap (\mathrm{Sing}(V_3)\cap P)=
\mathrm{Sing}(Q)\cap V_{2,2}.$$ Thus
$\mathrm{Sing}(Q)\cap\mathrm{Sing}(V_3)\neq\emptyset$ if $\dim
\mathrm{Sing}(Q)\geq 2$, i.e. if $n\geq 4$ in case (i), and $n\geq
3$ in case (ii). Finally, let us show that if
$\mathrm{Sing}(Q)\cap\mathrm{Sing}(V_3)\neq\emptyset$, then
$\mathrm{Sing}(X)\neq\emptyset$. Let
$x\in\mathrm{Sing}(Q)\cap\mathrm{Sing}(V_3)$. Then $x\in Q\cap
V_3=X\cup P$ and
$\mathbb{T}_x(Q)=\mathbb{T}_x(V_3)=\mathbb{P}_k^{n+2}$ since $\dim
Q=\dim V_3=n+1$. If $x\in P-X$, then $x$ is smooth for $Q\cap
V_3$, but $\mathbb{T}_x(Q\cap V_3)=\mathbb{T}_x(Q)\cap
\mathbb{T}_x(V_3)=\mathbb{P}_k^{n+2}$, a contradiction. Hence
$x\in X$. If $x$ would be a smooth point for $X$, then
$\mathbb{T}_xX\cong\mathbb{P}_k^n$. This implies that $x\in P\cap
X$. Moreover, since $\dim (P\cap X)=n-1$, we see that $\dim
(P\cap\mathbb{T}_xX)\geq n-1$. This shows that $\mathbb{T}_x(Q\cap
V_3)=\mathbb{T}_x(X\cup P)\neq\mathbb{P}_k^{n+2}$, a
contradiction. Thus $x\in\mathrm{Sing}(X)$. 

\medskip

\begin{rem}
All the surfaces in \cite[Thm. 4 (d)]{H1} are linked to a
$\mathbb{P}_{\mathbb{C}}^2$ by a complete intersection of type
$(3,2)$ as in case $(5)$ of the Theorem.
\end{rem}

\section{Some cases in any characteristic}

In this section, we extend to positive characteristic some special
cases of the Theorem stated in the Introduction. In particular,
let $k$ be an algebraically closed field of $\textrm{char}(k)\geq
0$.

First of all, let us note that in the case of normal surfaces we
can obtain the following

\begin{prop}\label{surface}
Let $S\subset\mathbb{P}_k^N$ be an irreducible reduced normal
surface defined over an algebraically closed field $k$. Denote by $H$ an ample and spanned line
bundle on $S$. If $\deg(S):=H^2=5$, then $(S,H)$ is one of the
following pairs:
\begin{enumerate}
\item $S$ is a covering $\pi:S\to\mathbb{P}_k^2$ of degree five
over $\mathbb{P}_k^2$ with
$H=\pi^*(\mathcal{O}_{\mathbb{P}_k^2}(1))$; \item $(S,H)$ is a
scroll over an elliptic curve, i.e. $S$ is a
$\mathbb{P}_k^1$-bundle over an elliptic curve $C$ and $H\cdot f
=1$, where $f$ is a fiber of the ruling $S\to C$; \item $(S,H)$ is
a scroll over $\mathbb{P}_k^1$, i.e. $S\cong\mathbb{P}(E)$ for an
ample vector bundle of rank two on $\mathbb{P}_k^1$ and $H=\xi$ is
the tautological line bundle; \item a cone over a smooth rational
curve of degree five; \item $(S,H)$ is a Gorenstein Del Pezzo
surface of degree five with $K_S\simeq -H$; \item
$S\subset\mathbb{P}_k^{4}$ is linked to a linear space
$P\cong\mathbb{P}_k^2$ by a complete intersection of type $(2,3)$,
i.e. $S\cup P=V_2\cap V_3$, where $V_3\subset\mathbb{P}_k^{4}$ is
a cubic hypersurface and $V_2\subset\mathbb{P}_k^{4}$ is a quadric
cone of rank $r$ with $3\leq r\leq 4$; moreover, if
$\mathrm{char}(k)\neq 2$, one of the following possibilities can
occurs:
\begin{enumerate}
\item[(i)] $r=3$, $S\simeq 5P_1$ and
$\mathrm{Cl}(V_2)\cong\mathbb{Z}[P_1]$; \item[(ii)] $r=4$,
$S\simeq 2P_1+3P_2$ and
$\mathrm{Cl}(V_2)\cong\mathbb{Z}[P_1]\oplus\mathbb{Z}[P_2]$,
\end{enumerate}
where $\mathrm{Cl}(V_2)$ is the class group of $V_2$ and the
$P_i$'s are linear spaces $\mathbb{P}_k^2$ such that $P_1=P$;
\item $S$ is a quintic hypersurface of $\mathbb{P}_k^{3}$.
\end{enumerate}
\end{prop}
\noindent\textit{Proof.} Since
$H^2=\deg\varphi(S)\cdot\deg\varphi$ is a prime number, where
$\varphi:S\to\mathbb{P}^N$ is the morphism associated to the
linear system $|H|$, we deduce that either (i) $\deg\varphi(S)=1,
\deg\varphi=5$, or (ii) $\deg\varphi(S)=5, \deg\varphi=1$. In case
(i) we see that $S$ is a covering $\pi:S\to\mathbb{P}^2$ of degree
five over $\mathbb{P}^2$ with
$H=\pi^*(\mathcal{O}_{\mathbb{P}^2}(1))$, while in case (ii) $H$
is very ample. From the $\Delta$-genus formula it follows that
\begin{equation}
4=\Delta + N - 2, \tag{$*'$}
\end{equation}
where $\Delta :=2+H^2-h^0(H)$ is the $\Delta$-genus of the
polarized pair $(S,H)$ and $h^0(H)=N+1$. Note that $N - 2\geq 1$
and $\Delta\geq 0$ by \cite[$\S$I (4.2)]{Fu1}. In particular, we
have $\Delta\leq 3$ with equality if and only if
$S\subset\mathbb{P}^{3}$ is a hypersurface of degree 5. So, from
$(*')$ we can assume that $\Delta\leq 2$. If $g(S,H)=0$, then by
\cite[$\S$I (3.4), (5.10) and (5.15)]{Fu1} we know that $(S,H)$ is
either
\begin{enumerate}
\item a smooth scroll over $\mathbb{P}_k^1$, i.e.
$S\cong\mathbb{P}(E)$ for an ample vector bundle $E$ on
$\mathbb{P}_k^1$ and $H$ is the tautological line bundle, or \item
a cone over a smooth rational curve of degree five.
\end{enumerate}
If $g(S,H)=1$ then by \cite[(1.2)(2)]{K} we know that $(S,H)$ is
one of the following two pairs:
\begin{enumerate}
\item[$(3)$] a scroll over an elliptic curve, i.e. $S$ is a
$\mathbb{P}^1$-bundle over an elliptic curve $C$ and $H\cdot f
=1$, where $f$ is a fiber of the ruling $S\to C$; or \item[$(4)$]
a Del Pezzo surface of degree five with $K_S\simeq -H$.
\end{enumerate}
Assume now that $g(S,H)\geq 2\geq\Delta$. From \cite[$\S$I
(3.5)]{Fu1} it follows that $g(S,H)=\Delta=2$ and by arguing as in
the proof of the Theorem we obtain the statement. \hfill $\square$

\begin{ex}
By using {\em MAGMA} \cite{magma}, the following commands 
allow us to construct examples of normal
varieties $X\subset\mathbb{P}_k^{n+2}$ linked to a linear space
$\mathbb{P}_k^n$ by a complete intersection of type $(2,3)$ as in
case $(5)$ of the Theorem.

\begin{verbatim}

//Fields K
K:=RationalField();
//or K:=GF(32003);

//Ambient space
P<[x]>:=ProjectiveSpace(K,6);
//or P<[x]>:=ProjectiveSpace(K,4);

//Random polynomials
f2:=Random(LinearSystem(P,2)); g2:=Random(LinearSystem(P,2));

//Function IsNormal()
function IsNormal(X);
 I:=Ideal(X);
 J:=Normalisation(I);
 P:=AmbientSpace(X);
 Y:=Scheme(P,(J[1][1]));
 test := Difference(X,Y) eq Difference(Y,X);
 if test then return true;
 else return false;
 end if;
end function;
\end{verbatim}

Case rank $3$:

\begin{verbatim}
Z:=Scheme(P,[x[1]*x[1]-x[2]*x[2]+x[3]*x[3],x[1]*f2-x[2]*f2+x[3]*g2]);
\end{verbatim}

Case rank $4$:

\begin{verbatim}
Z:=Scheme(P,[x[1]*x[1]-x[2]*x[2]+x[3]*x[3]-x[4]*x[4],x[1]*f2-x[2]*f2+x[3]*g2-x[4]*g2]);
\end{verbatim}

Therefore we can now define the variety $X$ as follows:

\begin{verbatim}
X:=PrimeComponents(Z)[1];
\end{verbatim}

Finally, we can check all the main properties of $X$ by putting

\begin{verbatim}
IsIrreducible(X); IsReduced(X); IsSingular(X); Dimension(X);
Degree(X); IsNormal(X);
\end{verbatim}

\end{ex}

\bigskip

On the other hand, when $n\geq 2$ and $X$ is a smooth variety, we
get the following

\begin{prop}\label{genus one}
Let $X$ be an irreducible and reduced smooth variety of dimension
$n\geq 2$ defined over an algebraically closed field $k$. Let $H$
be a very ample line bundle on $X$. If $g(X,H)=1$, then either
$\Delta (X,H)=1$ and $(X,H)$ is a Del Pezzo manifold, or $\Delta
(X,H)=n$ and $(X,H)\cong (\mathbb{P}(E),\xi)$, where $E$ is a
vector bundle over an elliptic curve and $\xi$ is the tautological
line bundle.
\end{prop}

\noindent\textit{Proof.} For $n=2$ the statement holds by
\cite[(1.2)(2)]{K}. Let $n\geq 3$. Applying the induction
hypothesis to $X_{n-1}\in |H|$, we see either (i) $\Delta
(X_{n-1},H_{n-1})=1$ and $(X_{n-1},H_{n-1})$ is a Del Pezzo
manifold, or (ii) $\Delta (X_{n-1},H_{n-1})=n-1$ and
$(X_{n-1},H_{n-1})$ is a scroll over an elliptic curve, where
$H_{n-1}:=H|_{X_{n-1}}$.

If $\Delta (X,H)\leq 1$ then by \cite[$\S$I (3.5)]{Fu1} we have
$\Delta (X,H)=g=1$. Thus from \cite[$\S$I (6.5)]{Fu1} it follows
that $(X,H)$ is a Del Pezzo manifold. So suppose that $\Delta
(X,H)\geq 2>1=g$.

Assume that $(X_{n-1},H_{n-1})$ is as in case (i). Then
$$h^0(H_{n-1}):=h^0(X_{n-1},H_{n-1})=n-1-\Delta (X_{n-1},H_{n-1})+H_{n-1}^{n-1}=n-2+H^n$$
and
$$h^0(H):=h^0(X,H)=n+H^n-\Delta (X,H)\leq n+H^n-2=h^0(H_{n-1}).$$ From the
exact sequence
$$0\to -H\to\mathcal{O}_X\to\mathcal{O}_{X_{n-1}}\to 0,$$ we see
$h^1(-H)=h^1(\mathcal{O}_X)$ since $h^1(\mathcal{O}_{X_{n-1}})=0$.
Moreover, for any integer $m\geq 1$ the exact sequence
$$0\to-(m+1)H\to -mH\to -mH_{X_{n-1}}\to 0$$ gives
$h^1(-(m+1)H)=...=h^1(-H)=h^1(\mathcal{O}_X)$ since
$h^1(-mH_{X_{n-1}})=0$ for any $m\geq 1$. By the Serre's duality
and the Serre vanishing theorem we conclude that
$h^1(\mathcal{O}_X)=0$. Therefore by using the following exact
sequence
$$0\to\mathcal{O}_X\to\mathcal{O}_X(H)\to\mathcal{O}_{X_{n-1}}(H)\to
0$$ we deduce that
$$h^0(H)=h^0(\mathcal{O}_X)+h^0(\mathcal{O}_{X_{n-1}}(H))=1+h^0(H_{n-1}),$$
but this leads to the contradiction $h^0(H)>h^0(H_{n-1})$.

So we can assume that $(X_{n-1},H_{n-1})$ is as in case (ii), i.e.
$(X_{n-1},H_{n-1})$ is a scroll over an elliptic curve $C$.
Consider a fiber $F\cong\mathbb{P}_k^{n-2}$ of $X_{n-1}\to C$.
Thus we have
$$(K_X+(n-1)H)|_F=[(K_X+(n-1)H)_{X_{n-1}}]_F\simeq (K_{X_{n-1}}+(\dim X_{n-1}-1)H_{n-1})_F=$$
$$=[K_{X_{n-1}}+(\dim X_{n-1})H_{n-1}]_F-[H_{X_{n-1}}]_F\simeq -[H_{X_{n-1}}]_F\simeq\mathcal{O}_{\mathbb{P}_k^{n-2}}(-1).$$
Hence $K_X+(n-1)H$ cannot be nef and by \cite[Theorem 1]{KS} we
see that $(X,H)$ is the projectivization of an ample vector bundle
of rank $n$ and $H=\xi$ is its tautological line bundle. Note that
the curve $X_1\in |H^{n-1}|$ is smooth and of genus one. Moreover,
$X_1\cdot F=H^{n-1}\cdot F=1$ and this gives $C\cong X_1$ and
$g(C)=g(X_1)=1$. Finally, by arguing as in the proof of the
Theorem, we can show that in this situation $\Delta(X,H)=n$.
\hfill $\square$

\bigskip

The above Proposition together with some results in \cite{Fu1},
allow us to obtain in the smooth case also the following

\begin{thm}
Let $X$ be an irreducible and reduced smooth variety of dimension
$n\geq 3$ defined over an algebraically closed field $k$.  
Let $H$ be an ample and spanned line
bundle on $X$. If $\deg(X):=H^n=5$, then $(X,H)$ is one of the
following pairs:
\begin{enumerate}
\item $X$ is a covering $\pi:X\to\mathbb{P}_k^n$ of degree five
over $\mathbb{P}_k^n$ with
$H=\pi^*(\mathcal{O}_{\mathbb{P}_k^n}(1))$; \item $(X,H)\cong
(\mathbb{P}(E),\xi)$ for some ample vector bundle $E$ on a smooth
curve $C$ with $g(C)\leq 1$, where $\xi$ is the tautological line
bundle on $\mathbb{P}(E)$; \item $n\leq 6$ and $X$ is a linear
section of the Grassmann variety parametrizing lines in
$\mathbb{P}_k^4$, embedded in $\mathbb{P}_k^9$ by Pl\"{u}cker
coordinates; \item $n=3$ and $X\subset\mathbb{P}_k^{5}$ is linked
to a linear space $P\cong\mathbb{P}_k^3$ by a complete
intersection of type $(2,3)$, i.e. $X\cup P=V_2\cap V_3$, where
$V_3\subset\mathbb{P}_k^{5}$ is a cubic hypersurface and
$V_2\subset\mathbb{P}_k^{5}$ is a quadric cone of
rank $4$; moreover, if $\mathrm{char}(k)\neq 2$ then $X\simeq 2P_1+3P_2$ is
a Weil divisor on $V_2$ and
$\mathrm{Cl}(V_2)\cong\mathbb{Z}[P_1]\oplus\mathbb{Z}[P_2]$,
where $\mathrm{Cl}(V_2)$ is the class group of $V_2$ and the
$P_i$'s are linear spaces $\mathbb{P}_k^3$ such that $P_1=P$;
\item $X$ is a quintic hypersurface of $\mathbb{P}_k^{n+1}$.
\end{enumerate}
\end{thm}

\noindent\textit{Proof.} By arguing as in Proposition
\ref{surface}, it follows that either $X$ is as in case (1) of the
statement, or $H$ is very ample on $X$. From the $\Delta$-genus
formula we have
\begin{equation}
4=\Delta + N - 2, \tag{$*''$}
\end{equation}
where $\Delta :=n+H^n-h^0(H)$ is the $\Delta$-genus of the
polarized pair $(X,H)$ and $h^0(H)=N+1$. Note that $N - 2\geq 1$
and $\Delta\geq 0$ by \cite[$\S$I (4.2)]{Fu1}. In particular, we
get $\Delta\leq 3$ with equality if and only if $X$ is as in case
(5) of the statement. So by $(*'')$ we can assume that $\Delta\leq
2$. If $g(S,H)\leq 1$, then by \cite[$\S$I (3.4), (5.10), (5.15),
(8.11) and (8.19)]{Fu1} and Proposition \ref{genus one}, we deduce
that $(X,H)$ is as in cases (2) and (3) of the statement. Finally,
assume that $g(S,H)\geq 2\geq\Delta$. By \cite[$\S$I (3.5)]{Fu1}
we obtain that $g(S,H)=\Delta=2$. Since $X$ is smooth, by arguing
as in the proof of the Theorem, 
it follows that $X$ is as in case (4) of the statement with $n=3$.
\hfill $\square$

\bigskip
\medskip

\noindent\textit{Acknowledgement}. I would like to thank Prof. A.
Laface for many interesting, stimulating and useful discussions
about algebraic geometry and for introducing me to the Magma
algebra system. This project was partially supported by Proyecto DIUC 2011, N.
211.013.036-1.0.

\bigskip
\medskip

Departamento de Matem\'{a}tica

Universidad de Concepci\'{o}n

Casilla 160-C, Concepci\'{o}n (Chile)

\smallskip

\textit{E-mail address:} atironi@udec.cl

\end{document}